\documentclass{article}
\usepackage{graphics}
\usepackage{graphicx}
\usepackage{latexsym}
\usepackage{amssymb}
\usepackage{amsmath}
\newtheorem{theorem}{Theorem}[section]
\newtheorem{lemma}[theorem]{Lemma}
\newtheorem{proposition}[theorem]{Proposition}

\newtheorem{preexample}{Example}[section]
\newenvironment{example}{\begin{preexample}}{\end{preexample}}
\newtheorem{preremark}{Remark}
\newenvironment{remark}{\begin{preremark}\rm}{\end{preremark}}

\newenvironment{proof}
  {{\bf Proof:}}
  {\qquad \hspace*{\fill} $\Box$}%

\newcommand{\fg}{\mathfrak{g}}%
\newcommand{\DC}{\mathcal{D}}%
\newcommand{\XC}{\mathcal{X}}%
\newcommand{\R}{\mathbb{R}}%

\begin{document}

\title{Controllability of Linear Control Systems on Solvable Lie Groups with Hyperbolic Drift}
\author{S. N. Stelmastchuk\\ Universidade Federal do Paran\'{a}\\Jandaia do Sul, Brazil, simnaos@ufpr.br}

\maketitle

\begin{abstract}
  We study the controllability of linear control systems on connected solvable Lie groups with hyperbolic drift. Under a suitable splitting assumption on the control directions, we show that the controllability problem can be reduced to the controllability of induced systems on the positive and negative hyperbolic components of the group. We then establish sufficient conditions ensuring controllability of these induced systems based on the ad-rank condition and suitable drift-reachability assumptions. As a consequence, we obtain a controllability criterion for the original system. The results provide a hyperbolic counterpart to existing controllability results for linear systems on solvable Lie groups and are illustrated by explicit examples.
\end{abstract}


\textbf{AMS 2010 subject classification}: 93C05, 93B05, 22E25

\section{Introduction}
Let $G$ be a connected solvable Lie group. Recall that a linear control system on $G$ is a control system of the form
\[
  \dot g(t)=\mathcal X(g(t))  +\sum_{i=1}^{m}u_i(t)X_i(g(t)), \tag{$\Sigma$}
\]
where $\mathcal X$ is a linear vector field, $X_1,\ldots,X_m$ are right-invariant vector fields, and $u=(u_1,\ldots,u_m)$ is an admissible control.

Linear control systems on Lie groups have been extensively studied during the last decades. Important contributions to the theory were obtained by Ayala et al. \cite{ayti}, Cardetti and Mittenhuber \cite{cardetti}, Jouan \cite{Jouan}, Da Silva \cite{daSilva}, among others. These works established several controllability, accessibility and structural results for this class of systems.

In the particular case of solvable Lie groups, Da Silva \cite{daSilva} obtained a controllability criterion under the assumption that the derivation associated with the drift has only eigenvalues with zero real part. The purpose of the present work is to investigate the complementary hyperbolic situation. More precisely, we assume that the derivation associated with the linear vector field $\mathcal X$ has no eigenvalues with zero real part.

Under this assumption, the Lie algebra $\mathfrak g$ admits the hyperbolic decomposition
\[
  \mathfrak g=\mathfrak g^+\oplus\mathfrak g^-,
\]
where $\mathfrak g^+$ and $\mathfrak g^-$ are the generalized eigenspaces associated with eigenvalues having positive and negative real parts, respectively. Let $G^+$ and $G^-$ denote the connected Lie subgroups associated with $\mathfrak g^+$ and $\mathfrak g^-$. Since $G$ is solvable and the drift is hyperbolic, we have the global factorization
\[
  G=G^+G^-=G^-G^+.
\]
Assuming that the control directions split independently according to this decomposition, the original system induces linear control systems $\Sigma^+$ and $\Sigma^-$ on $G^+$ and $G^-$, respectively.

Our first main result establishes that the controllability of the original system can be reduced to the controllability of the induced systems on the hyperbolic components. More precisely, we show that $\Sigma$ is controllable on $G$ whenever the induced systems $\Sigma^+$ and $\Sigma^-$ are controllable on $G^+$ and $G^-$, respectively (Theorem~\ref{thm:global_controllability_no_semigroup}).

The second part of the paper is devoted to obtaining controllability criteria for the induced systems. Our approach is based on the ad-rank condition. We prove that, under the split-control assumption, the ad-rank condition for the original system is equivalent to the corresponding ad-rank conditions on the hyperbolic components. Since the ad-rank condition alone does not imply controllability, we introduce the additional assumptions $(DR^-)$ and $(DR+^*)$. Under these hypotheses, we prove that the ad-rank condition implies controllability of the induced systems on $G^+$ and $G^-$ (Theorems~\ref{thm:controllability_Gplus_final} and \ref{thm:controllability_Gminus_final}). Combining these results with the factorization argument, we obtain our main controllability theorem (Theorem~\ref{thm:global_controllability_adrank_final}).

Finally, we present two examples illustrating the applicability of the theory. The first example is a hyperbolic linear control system on the Heisenberg group. The second example exhibits a genuinely nontrivial application of the ad-rank condition, where the control algebra does not coincide with one of the hyperbolic components, but the action of the derivation generates the missing directions.

The paper is organized as follows. In Section~2 we recall basic facts concerning linear vector fields and hyperbolic decompositions. In Section~3 we introduce the standing assumptions and establish the reduction of controllability from $G$ to the hyperbolic components. Section~4 contains the main controllability results based on the ad-rank condition. Finally, in Section~5 we present the examples illustrating the theory.

\section{Linear vector fields and hyperbolic decomposition}

Let \(G\) be a connected Lie group with Lie algebra \(\mathfrak g\). A vector field \(\mathcal X\) on \(G\) is called linear if its flow
\((\varphi_t)_{t\in\mathbb R}\) is a one-parameter subgroup of \(\operatorname{Aut}(G)\).

Associated with \(\mathcal X\) there is a derivation \(D\) of \(\mathfrak g\), defined by
\[
  D(Y)=-[\mathcal X,Y],\qquad Y\in\mathfrak g.
\]
We denote by
\[
  \mathfrak g^+,\qquad \mathfrak g^0,\qquad \mathfrak g^-
\]
the generalized eigenspaces of \(D\) associated, respectively, with eigenvalues having positive, zero and negative real parts. Thus
\[
  \mathfrak g=\mathfrak g^+\oplus\mathfrak g^0\oplus\mathfrak g^-.
\]
It is well known that the generalized eigenspaces of a derivation associated  with eigenvalues whose real parts have the same sign are Lie subalgebras. Therefore \(\mathfrak g^+\), \(\mathfrak g^0\), and \(\mathfrak g^-\) are Lie subalgebras of \(\mathfrak g\). Let \(G^+\), \(G^0\), and \(G^-\) be the connected Lie subgroups  associated with \(\mathfrak g^+\), \(\mathfrak g^0\), and \(\mathfrak g^-\), respectively.

A linear vector field \(\mathcal X\) is called {\em hyperbolic} if its associated derivation \(D\) has no eigenvalues with zero real part. 

In the hyperbolic case, \(\mathfrak g^0=\{0\}\). Hence the connected Lie subgroup associated with \(\mathfrak g^0\) is trivial, that is, \(G^0=\{e\}\). In consequence, $\mathfrak g=\mathfrak g^+\oplus\mathfrak g^-$.

When \(G\) is solvable and \(\mathcal X\) is hyperbolic, the global factorization [8, Proposition 2.9] is given by
\[
  G=G^+G^-=G^-G^+
\]
holds.

\begin{proposition}\label{invarianceofflow}
  Let $\XC$ be a linear vector field with flow $\varphi_t$. If $g \in G^{\pm}$, then $\varphi_t(g) \in G^{\pm}$ for all $t \in \R$. 
\end{proposition}
\begin{proof}
  Suppose that $g \in G^{\pm}$. Because $G^{\pm}$ is a connected Lie group there exists $X_1,\ldots, X_n \in \fg^{\pm}$ such that $g = \exp(X_1)\cdots \exp(X_n)$. Thus applying the flow $\varphi_t$ at equality we see that
	\[
	  \varphi_t(g) = \varphi_t( \exp(X_1))\cdots \varphi_t( \exp(X_n)) = \exp(e^{t\DC}X_1)\cdots \exp(e^{t\DC}X_n).
	\]
	Since $\fg^\pm$ is invariant by $\DC$, it follows that $e^{t\DC}X_1,\ldots e^{t\DC}X_n \in \fg^{\pm}$. In consequence, $\varphi_t(g) \in G^{\pm}$. 
\end{proof}

We will denote by $\varphi_t^\pm$ the restriction of $\varphi_t$ in $G^\pm$. It follows from Proposition 2.1 that the restrictions
\[
  \varphi_t^\pm:=\varphi_t|_{G^\pm}
\]
are well-defined. Since each \(\varphi_t\) is an automorphism of \(G\), each \(\varphi_t^\pm\) is an automorphism of \(G^\pm\). Therefore
\((\varphi_t^\pm)_{t\in\mathbb R}\) is a one-parameter subgroup of \(\operatorname{Aut}(G^\pm)\). Consequently, there exists a linear vector field \(\mathcal X^\pm\) on \(G^\pm\) whose flow is \((\varphi_t^\pm)_{t\in\mathbb R}\). The associated derivation is precisely the restriction
\[
  D^\pm=D|_{\mathfrak g^\pm}.
\]


\section{Controllability under hyperbolic decomposition and split controls}
\label{sec:hyperbolic_split_controls}

In this section we study the controllability of a linear control system on a connected solvable Lie group under a hyperbolic decomposition of the drift and an independent splitting of the control directions.

Standing assumption. Throughout the remainder of this paper, unless explicitly stated otherwise, \(G\) is a connected solvable Lie group and \(\mathcal X\) is a hyperbolic linear vector field on \(G\), with flow \((\varphi_t)_{t\in\mathbb R}\).
We assume that the associated hyperbolic decomposition satisfies
\[
    \mathfrak g=\mathfrak g^+\oplus\mathfrak g^-,
    \qquad
    G=G^+G^-=G^-G^+.
\]
Moreover, we assume that the control system has the split form
\[
    \dot g(t)
    =
    \mathcal X(g(t))
    +
    \sum_{\alpha=1}^{m_+}u_\alpha^+(t)X_\alpha^+(g(t))
    +
    \sum_{\beta=1}^{m_-}u_\beta^-(t)X_\beta^-(g(t)),
\]
where \(X_\alpha^+\) and \(X_\beta^-\) are right-invariant vector fields satisfying
\[
    X_\alpha^+(e)\in\mathfrak g^+,
    \qquad
    X_\beta^-(e)\in\mathfrak g^-.
\]
Finally, the admissible controls are assumed to split independently, in the sense that, whenever \(u^+\) and \(u^-\) are admissible for the corresponding families, the controls
\[
    (u^+,0),\qquad (0,u^-),\qquad (u^+,u^-)
\]
are admissible for the original system.

Since \(G^+\) and \(G^-\) are invariant under the flow \(\varphi_t\), and since the right-invariant fields generated by elements of \(\mathfrak g^\pm\) are tangent to \(G^\pm\), the system induces the following subsystems:
\begin{equation}\label{eq:condition_4}
    \Sigma^+:\quad
    \dot g(t)
    =
    \mathcal X^+(g(t))
    +
    \sum_{\alpha=1}^{m_+}u_\alpha^+(t)X_\alpha^+(g(t)),
    \qquad g(t)\in G^+,
\end{equation}
and
\begin{equation}\label{eq:condition_5}
    \Sigma^-:\quad
    \dot g(t)
    =
    \mathcal X^-(g(t))
    +
    \sum_{\beta=1}^{m_-}u_\beta^-(t)X_\beta^-(g(t)),
    \qquad g(t)\in G^-.
\end{equation}

We denote by \(\mathcal A(e)\) the reachable set from the identity for the original system, and by \(\mathcal A^*(e)\) the reachable set from the identity for the reverse system. Similarly, \(\mathcal A^\pm(e)\) and \(\mathcal A^{\pm,*}(e)\) denote the reachable sets from the identity for the induced systems and their reverse systems on \(G^\pm\).

Recall that, for linear control systems on Lie groups with right-invariant control vector fields, the control flow satisfies
\begin{equation}\label{eq:linear_flow_identity}
  \phi_{t,u}(g)=\phi_{t,u}(e)\varphi_t(g), \qquad g\in G,\ t\geq 0.
\end{equation}

\begin{theorem}\label{thm:global_controllability_no_semigroup}
  Under the standing assumptions, if the induced systems \(\Sigma^+\) and \(\Sigma^-\) are controllable on \(G^+\) and \(G^-\), respectively, then the original system \(\Sigma\) is controllable on \(G\).
\end{theorem}

\begin{proof}
  We first prove that every element of \(G\) is reachable from the identity for the original system. Let \(g\in G\). Since $G=G^-G^+$, there exist \(g^-\in G^-\) and \(g^+\in G^+\) such that
  \[
    g=g^-g^+.
  \]
  By controllability of \(\Sigma^-\) on \(G^-\), there exist \(s>0\) and an admissible control \(u^-\) such that
  \[
    \phi^-_{s,u^-}(e)=g^-.
  \]
  Since \(G^+\) is invariant under the linear flow \((\varphi)_{t\in\mathbb R}\), the element
  \[
    p:=\varphi^+_{-s}(g^+)
  \]
  belongs to \(G^+\). By controllability of \(\Sigma^+\) on \(G^+\), there exist \(r>0\) and an admissible control \(u^+\) such that
  \[
    \phi^+_{r,u^+}(e)=p.
  \]
  By the independent splitting of the controls, \((u^+,0)\) and \((0,u^-)\) are admissible controls for the original system. Applying first \((u^+,0)\), the  system reaches \(p\). Starting from \(p\), apply \((0,u^-)\) during time \(s\). Using the standard identity (\ref{eq:linear_flow_identity}) we obtain
  \[
    \phi_{s,(0,u^-)}(p) = \phi_{s,(0,u^-)}(e)\varphi_s(p).
  \]
  Since the trajectory generated by \((0,u^-)\) from the identity coincides with  the trajectory of \(\Sigma^-\), we have
  \[
    \phi_{s,(0,u^-)}(e)=\phi_{s,u^-}^-(e)=g^-.
  \]
  Moreover,
  \[
    \varphi_s(p)=\varphi_s(\varphi_{-s}^+(g^+))=g^+.
  \]
  Therefore
  \[
    \phi_{s,(0,u^-)}(p)=g^-g^+=g.
  \]
  Hence \(g\in\mathcal A(e)\). Since \(g\in G\) was arbitrary, $\mathcal A(e)=G$. It remains to prove that every element of \(G\) can be steered to the identity. Equivalently, we prove that
  \[
    \mathcal A^*(e)=G.
  \]
  The reverse system is again a linear control system, with linear vector field \(-\mathcal X\) and flow \(\varphi_{-t}\). Since \(G^+\) and \(G^-\) are invariant  under \(\varphi_t\), they are also invariant under the reverse flow. Moreover, because \(\Sigma^+\) and \(\Sigma^-\) are controllable, their reverse systems are also controllable. Applying the previous argument to the reverse system yields
  \[
    \mathcal A^*(e)=G.
  \]
  Thus every point of \(G\) is reachable from the identity and every point of \(G\) can be steered to the identity. Therefore \(\Sigma\) is controllable on \(G\).
 
\end{proof}
\begin{remark}
  The independent splitting of the controls is essential. It is not enough to know  that some projected systems on \(\mathfrak g^+\) and \(\mathfrak g^-\) are  controllable. The proof above uses the fact that trajectories of \(\Sigma^+\) and \(\Sigma^-\) are actual trajectories of the original system, obtained by setting one family of controls equal to zero.
\end{remark}

\begin{remark}
  The assumptions \(X_\alpha^+(e)\in\mathfrak g^+\) and \(X_\beta^-(e)\in\mathfrak g^-\) ensure that the corresponding right-invariant  fields are tangent to \(G^+\) and \(G^-\). Together with the invariance of these  subgroups under the drift flow, this makes the induced systems well-defined.
\end{remark}

\section{The ad-rank condition and controllability of the hyperbolic components}
\label{sec:adrank_hyperbolic_components}

In this section we work under the standing assumptions introduced in Section 3. Our aim is to show that, under suitable drift-reachability assumptions on the hyperbolic components, the ad-rank condition for the original system implies controllability of the induced systems on \(G^+\) and \(G^-\). Combined with Theorem~\ref{thm:global_controllability_no_semigroup}, this yields controllability of the original system on \(G\).

Let
\[
  h^+=\operatorname{Lie}\{X_1^+(e),\ldots,X_{m_+}^+(e)\},
  \qquad
  h^-=\operatorname{Lie}\{X_1^-(e),\ldots,X_{m_-}^-(e)\}.
\]
Under the standing assumptions, the full control algebra is
\[
  h=h^+\oplus h^-.
\]
Define
\[
    \mathcal D\mathfrak h:= \operatorname{span}\{\mathcal D^kY;\;Y\in\mathfrak h,\ k\geq0 \}.
\]
The system $\Sigma$ is said to satisfy the ad-rank condition on $G$ if
\[
    \mathcal D\mathfrak h=\mathfrak g.
\]
Similarly, for the positive subsystem, define
\[
    \mathcal D^+\mathfrak h^+ := \operatorname{span} \{(\mathcal D^+)^kY;\;Y\in\mathfrak h^+,\ k\geq0 \},
\]
where $\mathcal D^+:=\mathcal D|_{\mathfrak g^+}$. We say that $\Sigma^+$ satisfies the ad-rank condition on $G^+$ if
\[
    \mathcal D^+\mathfrak h^+=\mathfrak g^+.
\]
Analogously, for the negative subsystem, define
\[
    \mathcal D^-\mathfrak h^- := \operatorname{span} \{(\mathcal D^-)^kY;\;Y\in\mathfrak h^-,\ k\geq0 \},
\]
where $\mathcal D^-:=\mathcal D|_{\mathfrak g^-}$. We say that $\Sigma^-$ satisfies the ad-rank condition on $G^-$ if
\[
    \mathcal D^-\mathfrak h^-=\mathfrak g^-.
\]

\begin{proposition}\label{prop:adrank_decomposes_final}
  Assume that $\mathfrak h=\mathfrak h^+\oplus\mathfrak h^-$, $\mathfrak h^\pm\subset\mathfrak g^\pm$. Then $\mathcal D\mathfrak h
    = \mathcal D^+\mathfrak h^+ \oplus \mathcal D^-\mathfrak h^-$. Consequently, $\mathcal D\mathfrak h=\mathfrak g$ if and only if $ \mathcal D^+\mathfrak h^+=\mathfrak g^+$ and $\mathcal D^-\mathfrak h^-=\mathfrak g^-$.
\end{proposition}
\begin{proof}
  Since $\mathfrak g^+$ and $\mathfrak g^-$ are invariant under $\mathcal D$, we  have
  \[
    \mathcal D^k(\mathfrak g^+)\subset\mathfrak g^+,
    \qquad
    \mathcal D^k(\mathfrak g^-)\subset\mathfrak g^-,
  \]
  for every $k\geq0$.

  Let $Y\in\mathfrak h$. By assumption, there are unique elements $Y^+\in\mathfrak h^+$, $Y^-\in\mathfrak h^-$, such that $Y=Y^+ + Y^-$. Thus $ \mathcal D^kY = \mathcal D^kY^+ +\mathcal D^kY^-$. Moreover,
  \[
    \mathcal D^kY^+\in\mathfrak g^+,
    \qquad
    \mathcal D^kY^-\in\mathfrak g^-.
  \]
  Therefore $\mathcal D\mathfrak h = \mathcal D^+\mathfrak h^+ \mathcal D^-\mathfrak h^-$. Since $\mathfrak g=\mathfrak g^+\oplus\mathfrak g^-$, and $\mathcal D^+\mathfrak h^+\subset\mathfrak g^+$, $\mathcal D^-\mathfrak h^-\subset\mathfrak g^-$, the sum is direct. Hence
  \[
    \mathcal D\mathfrak h
    =
    \mathcal D^+\mathfrak h^+
    \oplus
    \mathcal D^-\mathfrak h^-.
  \]
  The final equivalence follows immediately from the decomposition $\mathfrak g=\mathfrak g^+\oplus\mathfrak g^-$.
\end{proof}

The previous proposition shows that, under the split-control assumption, the ad-rank condition on the original system is equivalent to the ad-rank conditions on the positive and negative hyperbolic components.

\begin{remark}  
  The previous proposition replaces the use of quotient projections. Since we do  not assume that $G^+$ or $G^-$ is normal in $G$, the maps that associate to an  element of $G$ its $G^+$ or $G^-$ component are not used as homomorphisms. The  relation between the ad-rank condition on $G$ and on the components is obtained directly from the vector space decomposition
  \[
    \mathfrak g=\mathfrak g^+\oplus\mathfrak g^-.
  \]
\end{remark}

We now prove that the ad-rank condition on each component implies controllability of the induced subsystem on that component. The proof uses two ingredients.

First, the ad-rank condition gives local openness properties of the reachable sets. More precisely, by the standard local accessibility results for linear control systems, if the ad-rank condition holds for the appropriate induced system, then the corresponding reachable sets have the required neighborhood properties around the identity.

Second, the hyperbolic behavior of the drift allows one to expand these local properties to the whole component.

\begin{lemma}\label{lem:adrank_openness_final}
The following assertions hold:
\begin{enumerate}
    \item If
    \[
        \mathcal D^+\mathfrak h^+=\mathfrak g^+,
    \]
    then $\mathcal A^+$ is open in $G^+$ and, for every $t>0$,
    $\mathcal A_t^{+,*}$ is a neighborhood of the identity in $G^+$.

    \item If
    \[
        \mathcal D^-\mathfrak h^-=\mathfrak g^-,
    \]
    then $\mathcal A^{-,*}$ is open in $G^-$ and, for every $t>0$,
    $\mathcal A_t^-$ is a neighborhood of the identity in $G^-$.
\end{enumerate}
\end{lemma}
\begin{proof}
  We prove the two assertions separately.

  \medskip
  \noindent
  \textit{Item 1.} Assume that $\mathcal D^+\mathfrak h^+=\mathfrak g^+$. This is precisely the ad-rank condition for the induced system $\Sigma^+$. By Proposition 6 in \cite{Jouan}, for all $t>0$ the reachable set $\mathcal{A}_t$ is a neighborhood of identity $e$. More precisely, in this setting, the ad-rank condition implies that $\mathcal A^+$ is open in $G^+$.

  Now consider the reverse system associated with $\Sigma^+$. Its linear vector field is $-\mathcal X|_{G^+}$ and the associated derivation is $\mathcal D^{+,*}=-\mathcal D^+$. Hence
  \[
    (\mathcal D^{+,*})^k=(-1)^k(\mathcal D^+)^k.
  \]
  Therefore
  \[
    \operatorname{span} \{ (\mathcal D^{+,*})^kY; \;Y\in\mathfrak h^+,\ k\geq0 \}
    =
    \operatorname{span} \{ (\mathcal D^+)^kY; \;Y\in\mathfrak h^+,\ k\geq0 \}.
  \]
  Since $\mathcal D^+\mathfrak h^+=\mathfrak g^+$, we also have $\mathcal D^{+,*}\mathfrak h^+=\mathfrak g^+$. Thus the reverse system on $G^+$ also satisfies the ad-rank condition.

  Applying Proposition 6 in \cite{Jouan} to the reverse system, we obtain, for every $t>0$, that $\mathcal A_t^{+,*}$, is a neighborhood of the identity in $G^+$. This proves the first assertion.

  \medskip
  \noindent
  \textit{Item 2.} Assume that \(D^-\mathfrak h^-=\mathfrak g^-\). This is the  ad-rank condition for the induced system \(\Sigma^-\) on \(G^-\). Hence, by Proposition 6 in \cite{Jouan}, for every \(t>0\), the reachable set \(\mathcal A_t^-\) is a neighborhood of the identity in \(G^-\). For the reverse system on \(G^-\), the associated derivation is $ D^{-,*}=-D^-$. Thus
  \[
    (D^{-,*})^k=(-1)^k(D^-)^k,
  \]
  and therefore $D^{-,*}\mathfrak h^- = D^-\mathfrak h^-=\mathfrak g^-$. Applying again Proposition 6 in \cite{Jouan} to the reverse system, we obtain that  \(\mathcal A^{-,*}\) is open in \(G^-\).
\end{proof}

\begin{remark}
  Lemma \ref{lem:adrank_openness_final} is the precise point where the ad-rank  condition is used. We do not assume directly that ad-rank implies global  controllability. Instead, we use only the local openness properties that follow from ad-rank and then combine them with the hyperbolic contraction-expansion  behavior on $G^+$ and $G^-$.
\end{remark}

We shall use the following drift-reachability assumptions:
\begin{equation*}\label{ass:dr+}
  \phi_t^-(g^-)\in\mathcal A^-(g^-),
  \qquad
  \forall g^-\in G^-,\ \forall t\in\mathbb R, \tag{DR-}
\end{equation*}
and
\[
  \phi_t^{+,*}(g^+)\in\mathcal A^{+,*}(g^+),
  \qquad
  \forall g^+\in G^+,\ \forall t\in\mathbb R. \tag{DR+*}
\]

Conditions $(DR^-)$ and $(DR^{+,*})$ are automatically satisfied whenever the control vector fields span the tangent space of the corresponding hyperbolic component. In this case, the drift can be compensated by the controls and every drift trajectory belongs to the reachable set from its initial point.

\begin{lemma}\label{lem:flow_reachability_backward}
  Let $\Sigma_H$ be a control system on a Lie group $H$ with drift flow $\psi_t$. Assume that, for every $h\in H$ and every $t\in\mathbb R$, $\psi_t(h)\in \mathcal A_H(h)$. Then, for every $s>0$, every $t>0$, and every $g\in H$,
  \[
    \psi_{-t}(\mathcal A_{H,s}(g))\subset \mathcal A_H(g).
  \]
\end{lemma}
\begin{proof}
  Let $h\in\mathcal A_{H,s}(g)$. Then there exists an admissible control $u$ such that $h=\phi_{s,u}(g)$. Thus $h\in\mathcal A_H(g)$. By hypothesis, applied to the point $h$ and to the time $-t$, we have
  \[
  \psi_{-t}(h)\in \mathcal A_H(h).
  \] 
  Concatenating a trajectory from $g$ to $h$ with a trajectory from $h$ to  $\psi_{-t}(h)$, we obtain $\psi_{-t}(h)\in\mathcal A_H(g)$. Therefore $\psi_{-t}(\mathcal A_{H,s}(g))\subset\mathcal A_H(g)$.
\end{proof}

We now prove the four reachability propositions needed to obtain controllability  on the components.

\begin{proposition}\label{prop:Gplus_Aplus_final}
  Assume that $G^+$ is nontrivial and that $\mathcal A^+$ is open in $G^+$. Then $ G^+\subset\mathcal A^+$.
\end{proposition}
\begin{proof}
  Let $g\in G^+$. Since $G^+$ is connected and nilpotent, the exponential map $\exp:\mathfrak g^+\to G^+$ is surjective. Hence there exists $X\in\mathfrak g^+$ such that $g=\exp X$.  Since every eigenvalue of $\mathcal D^+$ has positive real part, we have
  \[
    e^{-t\mathcal D^+}X\to0 \qquad\text{as }t\to+\infty.
  \]
  Therefore
  \[
    \varphi_{-t}^+(g) = \varphi_{-t}^+(\exp X)  = \exp(e^{-t\mathcal D^+}X)
  \]
  converges to the identity element of $G^+$ as $t\to+\infty$.
  Since $\mathcal A^+$ is open and contains the identity, there exists a  neighborhood $V$ of $e$ in $G^+$ such that $V\subset\mathcal A^+$. Thus, for $t>0$ sufficiently large, $\varphi_{-t}^+(g)\in\mathcal A^+$. By positive invariance of $\mathcal A^+$ under the flow $\varphi_t^+$, we obtain
  \[
    g = \varphi_t^+(\varphi_{-t}^+(g))  \in \varphi_t^+(\mathcal A^+) \subset \mathcal A^+.
  \]
  Therefore $G^+\subset\mathcal A^+$.
\end{proof}

\begin{proposition}\label{prop:Gminus_Aminus_final}
  Assume that $G^-$ is nontrivial and that, for every $s>0$, the set $\mathcal A_s^-$ is a neighborhood of the identity in $G^-$. Assume also that  for every $g\in G^-$ and every $t\in\mathbb R$, $\varphi_t^-(g)\in\mathcal A^-(g)$. Then $G^-\subset\mathcal A^-$.
\end{proposition}
\begin{proof}
  Let $g\in G^-$. Since $G^-$ is connected and nilpotent, the exponential map $\exp:\mathfrak g^-\to G^-$ is surjective. Hence there exists $X\in\mathfrak g^-$ such that $g=\exp X$. Since every eigenvalue of $\mathcal D^-$ has negative real part, we have
  \[
    e^{t\mathcal D^-}X\to0  \qquad\text{as }t\to+\infty.
  \]
  Therefore
  \[
    \varphi_t^-(g)  = \varphi_t^-(\exp X) = \exp(e^{t\mathcal D^-}X)
  \]
  converges to the identity element of $G^-$ as $t\to+\infty$. Choose $t>0$ sufficiently large and choose $s>t$ such that $\varphi_t^-(g)\in\mathcal A_s^-$. Using the negative-time reachability property with the zero control, we have
  \[
    g = \varphi_{-t}^-(\varphi_t^-(g))  \in \varphi_{-t}^-(\mathcal A_s^-)  \subset \mathcal A^-,
  \]
  where the last argument is given by Lemma~\ref{lem:flow_reachability_backward}. Therefore $G^-\subset\mathcal A^-$.
\end{proof}

\begin{proposition}\label{prop:Gplus_Aplus_reverse_flow_reachability}
  Assume that $G^+$ is nontrivial and that, for every $s>0$, the set  $\mathcal A_s^{+,*}$ is a neighborhood of the identity in $G^+$. Assume also  that for every $g\in G^+$ and every $t\in\mathbb R$,
  \[
    \varphi_t^{+,*}(g)\in\mathcal A^{+,*}(g).
  \]
  Then $G^+\subset\mathcal A^{+,*}$.
\end{proposition}
\begin{proof}
  For the reverse system on $G^+$, the associated derivation is
  \[
    \mathcal D^{+,*}=-\mathcal D^+.
  \]
  Thus every eigenvalue of $\mathcal D^{+,*}$ has negative real part. Let $g\in G^+$. Since $G^+$ is connected and nilpotent, there exists  $X\in\mathfrak g^+$ such that $g=\exp X$. Since
  \[
    e^{t\mathcal D^{+,*}}X\to0
    \qquad\text{as }t\to+\infty,
  \]
  we have $\varphi_t^{+,*}(g)\to e$. Choose $t>0$ sufficiently large. Since $\mathcal A_s^{+,*}$ is a neighborhood of the identity for every $s>0$, we may choose $s>0$ such that $\varphi_t^{+,*}(g)\in\mathcal A_s^{+,*}$. Using the negative-time reachability property with the zero control, we have
  \[
    g = \varphi_{-t}^{+,*}(\varphi_t^{+,*}(g)) \in \varphi_{-t}^{+,*}(\mathcal A_s^{+,*})  \subset \mathcal A^{+,*},
  \]
  where the last argument is given by Lemma~\ref{lem:flow_reachability_backward}. Therefore $G^+\subset\mathcal A^{+,*}$.
\end{proof}

\begin{proposition}\label{prop:Gminus_Aminus_reverse_final}
  Assume that $G^-$ is nontrivial and that $\mathcal A^{-,*}$ is open in $G^-$. Then
  \[
    G^-\subset\mathcal A^{-,*}.
  \]
\end{proposition}
\begin{proof}
  For the reverse system, the associated derivation is $\mathcal D^*=-\mathcal D$. Hence $\mathcal D^{-,*}:=\mathcal D^*|_{\mathfrak g^-}$ has only eigenvalues with positive real parts. 
  
  Let $g\in G^-$. Since $G^-$ is connected and nilpotent, there exists  $X\in\mathfrak g^-$ such that $g=\exp X$. Since the eigenvalues of $\mathcal D^{-,*}$ have positive real parts, $e^{-t\mathcal D^{-,*}}X\to\,0$ as $t\to+\infty$. Therefore
  \[
    \varphi_{-t}^{-,*}(g) = \exp(e^{-t\mathcal D^{-,*}}X)
  \]
  converges to the identity element of $G^-$. Since $\mathcal A^{-,*}$ is open and contains the identity, there exists a  neighborhood $V$ of $e$ in $G^-$ such that $V\subset\mathcal A^{-,*}$. Thus, for $t>0$ sufficiently large,
  \[
    \varphi_{-t}^{-,*}(g)\in\mathcal A^{-,*}.
  \]
  By positive invariance of $\mathcal A^{-,*}$ under the reverse flow, we obtain
  \[
    g = \varphi_t^{-,*}(\varphi_{-t}^{-,*}(g))  \in \varphi_t^{-,*}(\mathcal A^{-,*}) \subset \mathcal A^{-,*}.
  \]
  Therefore $G^-\subset\mathcal A^{-,*}$.
\end{proof}

We now prove that the ad-rank condition on each component implies controllability of the corresponding induced subsystem.

\begin{theorem}\label{thm:controllability_Gplus_final}
  Assume that \(G^+\) is nontrivial and that \((DR+*)\) holds. If \(D^+h^+=\mathfrak g^+\), then \(\Sigma^+\) is controllable on \(G^+\).
\end{theorem}
\begin{proof}
  Assume that $\mathcal D^+\mathfrak h^+=\mathfrak g^+$. By Lemma \ref{lem:adrank_openness_final}, the reachable set $\mathcal A^+$ is open in $G^+$ and, for every $t>0$, the reachable set $\mathcal A_t^{+,*}$ is a neighborhood of the identity in $G^+$.  By Proposition \ref{prop:Gplus_Aplus_final}, we have $G^+\subset\mathcal A^+$. Since $\mathcal A^+\subset G^+$, it follows that
  \[
    \mathcal A^+=G^+.
  \]
  By Proposition \ref{prop:Gplus_Aplus_reverse_flow_reachability}, we have $G^+\subset\mathcal A^{+,*}$. Again, since $\mathcal A^{+,*}\subset G^+$, it follows that
  \[
    \mathcal A^{+,*}=G^+.
  \]
  Therefore every element of $G^+$ is reachable from the identity and every element of $G^+$ can be steered to the identity. Hence $\Sigma^+$ is  controllable on $G^+$.
\end{proof}

\begin{theorem}\label{thm:controllability_Gminus_final}
  Assume that \(G^-\) is nontrivial and that \((DR-)\) holds. If  \(D^-h^-=\mathfrak g^-\), then \(\Sigma^-\) is controllable on \(G^-\).
\end{theorem}
\begin{proof}
  Assume that
  \[
    \mathcal D^-\mathfrak h^-=\mathfrak g^-.
  \]
  By Lemma \ref{lem:adrank_openness_final}, the reachable set $\mathcal A^{-,*}$  is open in $G^-$ and, for every $t>0$, the reachable set $\mathcal A_t^-$ is a  neighborhood of the identity in $G^-$.  By Proposition \ref{prop:Gminus_Aminus_final}, we have $G^-\subset\mathcal A^-$. Since $\mathcal A^-\subset G^-$, it follows that
  \[
    \mathcal A^-=G^-.
  \]
  By Proposition \ref{prop:Gminus_Aminus_reverse_final}, we have $G^-\subset\mathcal A^{-,*}$.  Since $\mathcal A^{-,*}\subset G^-$, it follows that
  \[
    \mathcal A^{-,*}=G^-.
  \]
  Therefore every element of $G^-$ is reachable from the identity and every element of $G^-$ can be steered to the identity. Hence $\Sigma^-$ is controllable on $G^-$.
\end{proof}

We now combine the previous results with the factorization argument.

\begin{theorem}\label{thm:global_controllability_adrank_final}
  Under the standing assumptions, assume that \((DR-)\) and \((DR+*)\) hold. If the original system satisfies the ad-rank condition
  \[
    Dh=\mathfrak g,
  \]
  then \(\Sigma\) is controllable on \(G\).
\end{theorem}
\begin{proof}
  By Proposition~\ref{prop:adrank_decomposes_final}, the condition $D\mathfrak h=\mathfrak g$ is equivalent to
  \[
    D^+\mathfrak h^+=\mathfrak g^+
    \qquad\text{and}\qquad
    D^-\mathfrak h^-=\mathfrak g^-.
  \]
  Therefore, by Theorems~\ref{thm:controllability_Gplus_final} and~\ref{thm:controllability_Gminus_final}, the induced systems \(\Sigma^+\) and \(\Sigma^-\) are controllable on \(G^+\) and \(G^-\),  respectively.

  Since the control directions split independently and
  \[
  G=G^+G^-=G^-G^+,
  \]
  Theorem~\ref{thm:global_controllability_no_semigroup} implies that the  original system \(\Sigma\) is controllable on \(G\).
\end{proof}

\begin{remark}
  The proof of Theorem~\ref{thm:global_controllability_adrank_final} has three distinct steps. First, the ad-rank condition on the full group decomposes into the ad-rank conditions on the hyperbolic components:
  \[
    D\mathfrak h=\mathfrak g
    \quad\Longleftrightarrow\quad
    D^+\mathfrak h^+=\mathfrak g^+
    \ \text{and}\
    D^-\mathfrak h^-=\mathfrak g^-.
  \]
  Second, these componentwise ad-rank conditions imply controllability of the induced systems on \(G^+\) and \(G^-\). Third, the global controllability of  \(\Sigma\) follows from Theorem~\ref{thm:global_controllability_no_semigroup},  using the global factorization
  \[
    G=G^+G^-=G^-G^+,
  \]
  the independent splitting of controls, and the standard flow identity for linear control systems with right-invariant control vector fields.
\end{remark}

\begin{remark}
  The theorem does not use projections from $G$ onto $G^+$ or $G^-$. This is  essential because, without assuming normality of $G^+$ or $G^-$, such component maps are not group homomorphisms in general. The argument relies only on the Lie algebra decomposition, the independent splitting of the controls, the global factorization of \(G\), and the standard flow identity for linear control systems.
\end{remark}

\section{Examples}

The purpose of this section is to illustrate how the controllability  criterion established in Theorem~4.10 can be applied in practice. We begin with a simple example where the ad-rank condition is verified directly. We then present a genuinely nontrivial situation in which the control directions do not generate the whole hyperbolic component, but the action of the derivation recovers the missing directions, showing the effectiveness of the ad-rank approach developed in the previous section.

\begin{example}\label{ex:heisenberg_hyperbolic}
  Let $\mathfrak g$ be the three-dimensional Heisenberg Lie algebra with basis
  \[
    \{X,Y,Z\}
  \]
  and nontrivial bracket $[X,Y]=Z$. Let $G$ be the connected simply connected Lie group associated with $\mathfrak g$. Since $\mathfrak g$ is nilpotent, $G$ is solvable.

  Define a linear map $\mathcal D:\mathfrak g\to\mathfrak g$ by
  \[
    \mathcal D X=X,\qquad
    \mathcal D Y=-2Y,\qquad
    \mathcal D Z=-Z.
  \]
  It is simple to show that $\mathcal D$ is a derivation. The eigenvalues of $\mathcal D$ are  $\{ 1,\,-2,\,-1\}$. Hence the associated linear vector field $\mathcal X$ is hyperbolic. Moreover,
  \[
    \mathfrak g^+=\mathbb R X, \qquad \mathfrak g^-=\operatorname{span}\{Y,Z\}.
  \]
  Let $G^+=\exp(\mathfrak g^+)$, $G^-=\exp(\mathfrak g^-)$. Since $G$ is connected, simply connected and nilpotent, the exponential map is  a diffeomorphism. Therefore every element of $G$ can be written as a product of an element of $G^+$ and an element of $G^-$, and also in the reverse order:
  \[
    G=G^+G^-=G^-G^+.
  \]
  Consider the linear control system
  \[
    \dot g(t) = \mathcal X(g(t))  + u_1(t)X(g(t)) + u_2(t)Y(g(t)) + u_3(t)Z(g(t)),
  \]
  where the control vector fields are right-invariant and $(u_1,u_2,u_3)\in\mathbb R^3$. The controls split independently as
  \[
    u^+(t)=u_1(t),
    \qquad
    u^-(t)=(u_2(t),u_3(t)).
  \]
  Moreover, $X\in\mathfrak g^+$, $Y,Z\in\mathfrak g^-$. Thus $\mathfrak h^+=\mathbb R X$, $\mathfrak h^-=\operatorname{span}\{Y,Z\}$,   and $\mathfrak h=\mathfrak h^+\oplus\mathfrak h^-$.
 
  We now verify the additional flow-reachability hypotheses required in Theorem \ref{thm:global_controllability_adrank_final}. First, consider the induced system on $G^-$. Since
  \[
    \mathfrak h^-=\operatorname{span}\{Y,Z\}=\mathfrak g^-,
  \]
  the right-invariant control vector fields $Y$ and $Z$ span the tangent space of $G^-$ at every point. Moreover, the control range is $\mathbb R^2$ on this  component. Let $g^-\in G^-$ and let $t\in\mathbb R$. If $t\geq0$, then $\varphi_t^-(g^-)$ is reached from $g^-$ by taking the zero control during time $t$. If $t<0$, consider the curve
  \[
    \gamma(\tau)=\varphi_{\tau t}^-(g^-),
    \qquad 0\leq \tau\leq1.
  \]
  Then  $\dot\gamma(\tau) = t\,\mathcal X^-(\gamma(\tau))$. On the other hand, the induced system on $G^-$ has the form
  \[
    \dot g(\tau)  = \mathcal X^-(g(\tau)) + u_2(\tau)Y(g(\tau)) + u_3(\tau)Z(g(\tau)).
  \]
  Since $Y$ and $Z$ span $T_{\gamma(\tau)}G^-$ for every $\tau$, and since the  controls are unrestricted, we can choose measurable controls $u_2,u_3$ such that
  \[
    u_2(\tau)Y(\gamma(\tau))  + u_3(\tau)Z(\gamma(\tau))  = (t-1)\mathcal X^-(\gamma(\tau)).
  \]
  Therefore $\gamma$ is an admissible trajectory of $\Sigma^-$ from $g^-$ to $\varphi_t^-(g^-)$. Hence 
  \[
    \varphi_t^-(g^-)\in\mathcal A^-(g^-), \qquad \forall g^-\in G^-, \quad \forall t\in\mathbb R.
  \]
  
  Similarly, for the reverse system on $G^+$, we have
  \[
    \mathfrak h^+=\mathbb RX=\mathfrak g^+.
  \]
  Thus the right-invariant control vector field $X$ spans the tangent space of  $G^+$ at every point, and the control range is $\mathbb R$. Repeating the same  argument for the reverse drift, we obtain
  \[
    \varphi_t^{+,*}(g^+)\in\mathcal A^{+,*}(g^+),
    \qquad
    \forall g^+\in G^+,
    \quad
    \forall t\in\mathbb R.
  \]
  Therefore the additional hypotheses of Theorem \ref{thm:global_controllability_adrank_final} are satisfied.
  
  We now verify the ad-rank condition. Since $\mathcal D X=X$, we have  
  \[
    \mathcal D^+\mathfrak h^+
    =
    \operatorname{span}\{X\}
    =
    \mathfrak g^+.
  \]
  Similarly, $\mathcal D Y=-2Y$, $\mathcal D Z=-Z$, and therefore
  \[
    \mathcal D^-\mathfrak h^-
    =
    \operatorname{span}\{Y,Z\}
    =
    \mathfrak g^-.
  \]
  Consequently, $\mathcal D\mathfrak h = \mathfrak g$. Hence the system satisfies the ad-rank condition on $G$.

  By Theorems \ref{thm:controllability_Gplus_final} and \ref{thm:controllability_Gminus_final}, the induced systems on $G^+$ and  $G^-$ are controllable. Since the controls split independently and
  \[
    G=G^+G^-=G^-G^+,
  \]  
  Theorem \ref{thm:global_controllability_adrank_final} implies that the original system is controllable on $G$.
\end{example}

This example illustrates a genuinely nontrivial application of the  ad-rank condition. In contrast with the previous examples, the control  algebra does not coincide with the whole negative component.  Nevertheless, the action of the derivation generates the missing
directions.

\begin{example}\label{ex:fivedimensionalalgebra}
  Let \(\mathfrak g\) be the five-dimensional real Lie algebra with basis
  \[
    \{X_1,X_2,X_3,Y_1,Y_2\}
  \]
  and only nonzero bracket
  \[
    [X_1,X_2]=X_3.
  \]
  Thus \(\operatorname{span}\{X_1,X_2,X_3\}\) is a three-dimensional Heisenberg Lie algebra and \(\operatorname{span}\{Y_1,Y_2\}\) is an abelian ideal. In  particular, \(\mathfrak g\) is nilpotent and hence solvable.

  Let \(G\) be the connected and simply connected Lie group associated with \(\mathfrak g\). Define a linear map \(D:\mathfrak g\to\mathfrak g\) by
  \[
    DX_1=X_1,\qquad
    DX_2=2X_2,\qquad
    DX_3=3X_3,
  \]
  and
  \[
    DY_1=-Y_1+Y_2,\qquad
    DY_2=-2Y_2.
  \]
  Then \(D\) is a derivation. Indeed, the only nonzero bracket to be checked is \([X_1,X_2]=X_3\), and
  \[
    D[X_1,X_2]=DX_3=3X_3,
  \]
  while
  \[
    [DX_1,X_2]+[X_1,DX_2]
    =
    [X_1,X_2]+[X_1,2X_2]
    =
    3X_3.
  \]
  All brackets involving \(Y_1\) and \(Y_2\) are zero, so the derivation property is satisfied. The eigenvalues of \(D\) are $\{1,\ 2,\ 3,\ -1,\ -2\}$. Hence the associated linear vector field \(\mathcal X\) is hyperbolic. Moreover,
  \[
    \mathfrak g^+ = \operatorname{span}\{X_1,X_2,X_3\},
    \qquad
    \mathfrak g^- = \operatorname{span}\{Y_1,Y_2\}.
  \]
  Let $G^+=\exp(\mathfrak g^+)$, $G^-=\exp(\mathfrak g^-)$. Since \(G\) is connected, simply connected and nilpotent, the exponential map is a diffeomorphism. Moreover, \(\mathfrak g=\mathfrak g^+\oplus\mathfrak g^-\) with \([\mathfrak g^+,\mathfrak g^-]=0\). Therefore
  \[
    G=G^+G^-=G^-G^+.
  \]
  Consider the linear control system
  \[
    \dot g(t)
    =
    \mathcal X(g(t))
    +
    u_1(t)X_1(g(t))
    +
    u_2(t)X_2(g(t))
    +
    u_3(t)X_3(g(t))
    +
    v(t)Y_1(g(t)),
  \]
  where the control vector fields are right-invariant and the control range is  \(\mathbb R^4\). The controls split independently as
  \[
    u^+(t)=(u_1(t),u_2(t),u_3(t)),
    \qquad
    u^-(t)=v(t).
  \]  
  Thus
  \[
    h^+
    =
    \operatorname{span}\{X_1,X_2,X_3\}
    =
    \mathfrak g^+, 
    \mbox{ and }
    h^-=\operatorname{span}\{Y_1\}
    \subsetneq
    \mathfrak g^-.
  \]
  We now verify the ad-rank condition on the negative component. Since
  \[
    D^-Y_1=-Y_1+Y_2,
  \]
  we have
  \[
    D^-h^-
    =
    \operatorname{span}\{Y_1,D^-Y_1\}
    =
    \operatorname{span}\{Y_1,-Y_1+Y_2\}
    =
    \operatorname{span}\{Y_1,Y_2\}
    =
    \mathfrak g^-.
  \]
  Therefore the ad-rank condition on \(G^-\) holds although $h^-\neq \mathfrak g^-$.
  
  On the positive component, since \(h^+=\mathfrak g^+\), we immediately have
  \[
    D^+h^+=\mathfrak g^+.
  \]
  Consequently,
  \[
    Dh
    =
    D^+h^+\oplus D^-h^-
    =
    \mathfrak g^+\oplus\mathfrak g^-
    =
    \mathfrak g.
  \]
  Thus the full ad-rank condition holds.

  It remains to verify the additional drift-reachability hypotheses of  Theorem~\ref{thm:global_controllability_adrank_final}. For the positive component,  we have \(h^+=\mathfrak g^+\). Hence the right-invariant control vector fields  \(X_1,X_2,X_3\) span the tangent space of \(G^+\) at every point. Since the controls are unrestricted, the same compensation argument used in the previous  examples shows that
  \[
    \phi_t^{+,*}(g^+)\in \mathcal A^{+,*}(g^+),
    \qquad
    \forall g^+\in G^+,\ \forall t\in\mathbb R.
  \]

  For the negative component, \(G^-\) is abelian and may be identified with \(\mathbb R^2\). In the basis \(\{Y_1,Y_2\}\), the induced negative system is the classical linear control system
  \[
    \dot x(t)=D^-x(t)+v(t)Y_1.
  \]
  The Kalman controllability matrix is
  \[
    \bigl[Y_1,\ D^-Y_1\bigr].
  \]
  Since
  \[
    \operatorname{span}\{Y_1,D^-Y_1\}
    =
    \operatorname{span}\{Y_1,Y_2\}
    =
    \mathfrak g^-,
  \]
  this linear system is controllable on \(G^-\simeq\mathbb R^2\). In particular,  for every \(g^-\in G^-\) and every \(t\in\mathbb R\),
  \[
    \phi_t^-(g^-)\in \mathcal A^-(g^-).
  \]
  Therefore all hypotheses of Theorem \ref{thm:global_controllability_adrank_final} are satisfied. Hence the original system is controllable on \(G\).
\end{example}

\end{document}